\documentclass[12pt]{article}%
\usepackage{graphicx}
\usepackage[intlimits]{amsmath}
\usepackage{latexsym}
\usepackage{amsfonts}
\usepackage{amssymb}%
 \makeatletter
\newfont{\footsc}{cmcsc10 at 8truept}
\newfont{\footbf}{cmbx10 at 8truept}
\newfont{\footrm}{cmr10 at 10truept}
\newtheorem{theorem}{Theorem}

\newtheorem{lemma}[theorem]{Lemma}
\newtheorem{fact}[theorem]{Fact}

\newenvironment{proof}[1][Proof]{\noindent{\textbf {#1}  }}  {\hfill$\Box$\bigskip}

\begin{document}

\title{A spectral condition for odd cycles in graphs}
\author{Vladimir Nikiforov\\{\small Department of Mathematical Sciences, University of Memphis, Memphis TN
38152}\\{\small email: vnikifrv@memphis.edu}}
\maketitle

\begin{abstract}
Let $G$ be a graph of sufficiently large order $n,$ and let the largest
eigenvalue $\mu\left(  G\right)  $ of its adjacency matrix satisfies
$\mu\left(  G\right)  >\sqrt{\left\lfloor n^{2}/4\right\rfloor }.$ Then $G$
contains a cycle of length $t$ for every $t\leq n/320.$

This condition is sharp: the complete bipartite graph $T_{2}\left(  n\right)
$ with parts of size $\left\lfloor n/2\right\rfloor $ and $\left\lceil
n/2\right\rceil $ contains no odd cycles and its largest eigenvalue is equal
to $\sqrt{\left\lfloor n^{2}/4\right\rfloor }.$

This condition is stable: if $\mu\left(  G\right)  $ is close to
$\sqrt{\left\lfloor n^{2}/4\right\rfloor }$ and $G$ fails to contain a cycle
of length $t$ for some $t\leq n/321,$ then $G$ resembles $T_{2}\left(
n\right)  .$\medskip



\end{abstract}

\section*{Introduction}

This note is part of an ongoing project aiming to build extremal graph theory
on spectral grounds, see, e.g., \cite{BoNi07} and \cite{Nik02,Nik07h}.

It is known (\cite{Nik07d}, \cite{Nos70}) that if $G$ is a graph of order $n$
and the largest eigenvalue $\mu\left(  G\right)  $ of its adjacency matrix
satisfies $\mu\left(  G\right)  >\sqrt{\left\lfloor n^{2}/4\right\rfloor },$
then a triangle exists in $G$.

Here we show that the same premises imply the existence of other cycles as well.

\begin{theorem}
\label{th1}Let $G$ be a graph of sufficiently large order $n$ with $\mu\left(
G\right)  >\sqrt{\left\lfloor n^{2}/4\right\rfloor }.$ Then $G$ contains a
cycle of length $t$ for every $t\leq n/320.$
\end{theorem}

Write $T_{2}\left(  n\right)  $ for the complete bipartite graph with parts of
size $\left\lfloor n/2\right\rfloor $ and $\left\lceil n/2\right\rceil .$ Note
that $T_{2}\left(  n\right)  $ contains no odd cycles and $\mu\left(
T_{2}\left(  n\right)  \right)  =\sqrt{\left\lfloor n^{2}/4\right\rfloor };$
thus, Theorem \ref{th1} gives a sharp spectral condition for the existence of
odd cycles.

Moreover, there is stability in this condition: if $\mu\left(  G\right)  $ is
close to $\sqrt{\left\lfloor n^{2}/4\right\rfloor }$ and $G$ fails to contain
a cycle of length $t$ for some $t\leq n/321,$ then $G$ resembles $T_{2}\left(
n\right)  .$ Here is a precise form of this statement.

\begin{theorem}
\label{th2}Let $0<\theta<2^{-16}$ and $n$ be sufficiently large. For every
graph $G$ of order $n$ with $\mu\left(  G\right)  >\left(  1/2-\theta\right)
n,$ one of the following conditions holds;

(i) $G$ contains a cycle of length $t$ for every $t\leq n/321;$

(ii) there exists an induced bipartite subgraph $G_{0}\subset G$ satisfying
$\left\vert G_{0}\right\vert >\left(  1-4\theta^{1/3}\right)  n$ and
$\delta\left(  G_{0}\right)  >\left(  1/2-7\theta^{1/3}\right)  n.$
\end{theorem}

The proofs of Theorems \ref{th1} and \ref{th2} are based on three results of
independent interest.

\begin{lemma}
\label{lev1}Let $G$ be a graph of order $n$ with minimum degree $\delta$ and
$\mu\left(  G\right)  =\mu.$ If $\left(  x_{1},\ldots,x_{n}\right)  $ is a
unit eigenvector to $\mu,$ then%
\[
\min\left\{  x_{1},\ldots,x_{n}\right\}  \leq\sqrt{\frac{\delta}{\mu
^{2}+\delta n-\delta^{2}}}.
\]

\end{lemma}

\begin{lemma}
\label{lev2}Let $G$ be a graph of order $n$ with $\mu\left(  G\right)  =\mu.$
If $\left(  x_{1},\ldots,x_{n}\right)  $ is a unit eigenvector to $\mu$ and
$u$ is\ a vertex satisfying $x_{u}=\min\left\{  x_{1},\ldots,x_{n}\right\}  ,$
then%
\[
\frac{\mu\left(  G-u\right)  }{n-1}>\frac{\mu\left(  G\right)  }{n}\left(
1+\frac{1}{n-1}\left(  1-nx_{u}^{2}-\frac{1}{n-1}\right)  \right)  .
\]

\end{lemma}

Combining these two lemmas, we get Theorem \ref{thv3} below. We hope that this
technical statement can be used in other spectral extremal problems.

\begin{theorem}
\label{thv3}Let $0<4\alpha\leq1$, $0<2\beta\leq1,$ $1/2-\alpha/4\leq\gamma<1,$
$K\geq0,$ and $n\geq\left(  42K+4\right)  /\alpha^{2}\beta$. If $G$ is a graph
of order $n$ with
\[
\mu\left(  G\right)  >\gamma n-K/n\text{ \ \ and \ \ }\delta\left(  G\right)
\leq\left(  \gamma-\alpha\right)  n,
\]
then there exists an induced subgraph $H\subset G\ $with $\left\vert
H\right\vert \geq\left(  1-\beta\right)  n,$ satisfying one of the following conditions:

(i) $\mu\left(  H\right)  >\gamma\left(  1+\beta\alpha/2\right)  \left\vert
H\right\vert ;$

(ii) $\mu\left(  H\right)  >\gamma\left\vert H\right\vert $ and $\delta\left(
H\right)  >\left(  \gamma-\alpha\right)  \left\vert H\right\vert .$
\end{theorem}

\section*{Proofs}

We start with some notation and results needed for our proofs.

Our graph-theoretical notation follows \cite{Bol98}. Specifically, given a
graph $G,$ we write:

- $\left\vert G\right\vert $ for the number of vertices of $G;$

- $E\left(  G\right)  $ for the edge set of $G;$

- $k_{3}\left(  G\right)  $ for the number of triangles of $G;$

- $d\left(  u\right)  $ for the degree of a vertex $u;$

- $\Gamma\left(  u\right)  $ for the set of neighbors of a vertex $u;$

- $\delta\left(  G\right)  $ for the minimum degree of $G.$

The following fact is a reduced version of Theorem 1 of \cite{NiSc06}.

\begin{fact}
\label{tNS} Let $G$ be a nonbipartite graph of sufficiently large order $n,$
and let $\delta\left(  G\right)  \geq n/3.$ Then $C_{t}\subset G$ for every
integer $t\in\left[  4,\delta\left(  G\right)  +1\right]  .\hfill\square$
\end{fact}

The following facts are particular cases of Theorems 2 and 4 in \cite{BoNi07}.

\begin{fact}
\label{leNSMM}If $G$ is a graph of order $n,$ then $k_{3}\left(  G\right)
\geq\left(  \mu\left(  G\right)  /n-1/2\right)  n^{3}/12.\hfill\square$
\end{fact}

\begin{fact}
\label{tstab}Let $0<\theta\leq2^{-16}$ and let $G$ be a triangle-free graph of
order $n$ with $\mu\left(  G\right)  \geq\left(  1/2-\theta\right)  n.$ Then
there exists an induced bipartite graph $H\subset G$ satisfying $\left\vert
H\right\vert >\left(  1-3\theta^{1/3}\right)  n$ and $\delta\left(  H\right)
>\left(  1/2-6\theta^{1/3}\right)  n.\hfill\square$
\end{fact}

\bigskip

\subsubsection*{\textbf{Proof of Lemma \ref{lev1}}}

\begin{proof}
[\textbf{ }]Set $\sigma=\min\left\{  x_{1},\ldots,x_{n}\right\}  .$ If
$\sigma=0,$ the assertion holds trivially, so we assume that $\sigma>0.$ This
implies also that $\delta>0$. Taking $u\in V\left(  G\right)  $ to satisfy
$d\left(  u\right)  =\delta,$ we have
\begin{align*}
\mu^{2}\sigma^{2}  &  \leq\mu^{2}x_{u}^{2}=\left(
{\textstyle\sum\limits_{i\in\Gamma\left(  u\right)  }}
x_{i}\right)  ^{2}\leq\delta%
{\textstyle\sum\limits_{i\in\Gamma\left(  u\right)  }}
x_{i}^{2}\leq\delta\left(  1-%
{\textstyle\sum\limits_{i\in V\left(  G\right)  \backslash\Gamma\left(
u\right)  }}
x_{i}^{2}\right) \\
&  \leq\delta\left(  1-\left(  n-\delta\right)  \sigma^{2}\right)
=\delta-\left(  \delta n-\delta^{2}\right)  \sigma^{2},
\end{align*}
implying that $\left(  \mu^{2}+\delta n-\delta^{2}\right)  \sigma^{2}%
\leq\delta,$ and the desired inequality follows.
\end{proof}

\subsubsection*{\textbf{Proof of Lemma \ref{lev2}}}

\begin{proof}
[\textbf{ }]Set for short $c=1-nx_{u}^{2}$ and $\mu=\mu\left(  G\right)  .$ We
have%
\[
\mu x_{u}=\sum_{v\in\Gamma\left(  u\right)  }x_{v}\text{ \ \ and \ \ }%
\mu=2\sum_{vw\in E\left(  G\right)  }x_{v}x_{w}.
\]
Hence, by Rayleigh's principle, we obtain
\[
\mu=2\sum_{vw\in E\left(  G-u\right)  }x_{v}x_{w}+2x_{u}\sum_{v\in
\Gamma\left(  u\right)  }x_{v}\leq\mu\left(  G-u\right)  \left(
1-x^{2}\right)  +2x_{u}^{2}\mu\left(  G\right)  ,
\]
implying that%
\begin{equation}
\frac{\mu\left(  G-u\right)  }{n-1}\geq\frac{\mu\left(  G\right)  }{n-1}%
\cdot\frac{1-2x_{u}^{2}}{1-x_{u}^{2}}=\frac{\mu\left(  G\right)  }{n-1}%
\cdot\left(  \frac{n-2+2c}{n-1+c}\right)  . \label{in1}%
\end{equation}
On the other hand, in view of $0\leq c\leq1,$ we find that%
\begin{align*}
n\left(  \frac{n-2+2c}{n-1+c}\right)  -n+1-c+\frac{1}{n-1}  &  =\frac
{-1+c+cn}{n-1+c}-c+\frac{1}{n-1}\\
&  =-\frac{\left(  1-c\right)  ^{2}}{n-1+c}+\frac{1}{n-1}\geq0.
\end{align*}
Hence, inequality (\ref{in1}) implies that
\[
\frac{\mu\left(  G-u\right)  }{n-1}\geq\frac{\mu\left(  G\right)  }{n-1}%
\cdot\left(  \frac{n-2+2c}{n-1+c}\right)  \geq\frac{\mu\left(  G\right)  }%
{n}\cdot\left(  1+\frac{c}{n-1}-\frac{1}{\left(  n-1\right)  ^{2}}\right)  ,
\]
completing the proof.
\end{proof}

\subsubsection*{Proof of Theorem \textbf{\ref{thv3}}}

\begin{proof}
[\textbf{ }]Let $\alpha,\beta,\gamma,K,n,$ and the graph $G$ satisfy the
conditions of the theorem. We immediately see that
\[
n\geq\frac{42K+4}{\alpha^{2}\beta}>\max\left\{  \frac{15}{\left(
1-\beta\right)  \alpha},\frac{1}{\beta},\frac{84K}{\alpha^{2}}\right\}  .
\]

Define a sequence of graphs $G_{0},\ldots,G_{k}$ by the following procedure
$\mathcal{P}$:

\qquad\textbf{begin}

\qquad\ \ \ \emph{set }$G_{0}=G;$

$\qquad\ \ \ \emph{set}$ $k=0;$

\qquad\ \ \ \textbf{while }$\delta\left(  G_{i}\right)  \leq\left(
\gamma-\alpha\right)  \left(  n-k\right)  $ $\emph{and}$ $k<\left\lfloor \beta
n\right\rfloor $ \textbf{do }

\qquad\ \ \ \textbf{begin}

\qquad\ \ \ \ \ \ \ \emph{select a unit eigenvector }$\left(  x_{1}%
,\ldots,x_{n-k}\right)  $ \emph{to }$\mu\left(  G_{k}\right)  ;$

\qquad\ \ \ \ \ \ \ \emph{select a vertex }$u_{k}\in V\left(  G_{k}\right)  $
\emph{such that }$x_{u_{k}}=\min\left\{  x_{1},\ldots,x_{n-k}\right\}  ;$

\qquad\ \ \ \ \ \ \ \emph{set }$G_{k+1}=G_{k}-u_{k};$

\qquad\ \ \ \ \ \ \ \emph{add }$1$\emph{ to }$k;$

\qquad\ \ \ \textbf{end};

\qquad\textbf{end}.

Let $H=G_{k}$ and note that
\[
\left\vert H\right\vert =n-k\geq n-\left\lfloor \beta n\right\rfloor
\geq\left(  1-\beta\right)  n.
\]
We shall show that
\begin{equation}
\mu\left(  H\right)  >\gamma\left(  1+\frac{4k\alpha}{7n}\right)  \left\vert
H\right\vert . \label{bnd}%
\end{equation}
To this end, we first prove by induction on $i$ that
\begin{equation}
\frac{\mu\left(  G_{i}\right)  }{n-i}\geq\left(  1+\frac{3i\alpha}{5n}\right)
\frac{\mu\left(  G\right)  }{n} \label{in}%
\end{equation}
for every $i=0,\ldots,k.$

The assertion is trivially true for $i=0.$ Let $0\leq i\leq k-1$ and assume
that (\ref{in}) holds for $i;$ we shall prove that it also holds for $i+1$.
Set $\delta=\delta\left(  G_{i}\right)  $, $\mu=\mu\left(  G_{i}\right)  ,$
and note first that
\begin{align}
\delta &  \leq\left(  \gamma-\alpha\right)  \left(  n-i\right)  ,\label{mind}%
\\
\mu &  \geq\left(  n-i\right)  \left(  1+\frac{3i\alpha}{5n}\right)  \frac
{\mu\left(  G\right)  }{n}>\left(  n-i\right)  \left(  \gamma-\frac{K}{n^{2}%
}\right)  . \label{minm}%
\end{align}
Let $\left(  x_{1},\ldots,x_{n-i}\right)  $ be a unit eigenvector to\emph{
}$\mu,$ and let $u\in V\left(  G_{i}\right)  $ satisfy $x_{u}=\min\left\{
x_{1},\ldots,x_{n-i}\right\}  .$ Then Lemma \ref{lev1} implies that%
\[
x_{u}^{2}\leq\frac{\delta}{\mu^{2}+\left(  n-i\right)  \delta-\delta^{2}}.
\]
Noting that the right-hand side increases with $\delta$ and decreases with
$\mu,$ in view of (\ref{mind}) and (\ref{minm}), we find that%
\begin{align*}
x_{u}\left(  n-i\right)   &  \leq\frac{\left(  \gamma-\alpha\right)  \left(
n-i\right)  ^{2}}{\left(  n-i\right)  ^{2}\left(  \gamma-K/n^{2}\right)
^{2}+\left(  n-i\right)  \left(  \gamma-\alpha\right)  \left(  n-i\right)
-\left(  \gamma-\alpha\right)  ^{2}\left(  n-i\right)  ^{2}}\\
&  <\frac{\gamma-\alpha}{\gamma^{2}-2\gamma K/n^{2}+\gamma-\alpha
+2\gamma\alpha-\alpha^{2}}<\frac{\gamma-\alpha}{\gamma\left(  \gamma
+2\alpha\right)  -\alpha-2K/n^{2}+\gamma-\alpha^{2}}\\
&  \leq\frac{\gamma-\alpha}{\left(  1/2-\alpha/4\right)  \left(
1/2-\alpha/4+2\alpha\right)  -\alpha-2K/n^{2}+\gamma-\alpha^{2}}\\
&  <\frac{\gamma-\alpha}{-2K/n^{2}+\gamma-\alpha^{2}}<\frac{\gamma-\alpha
}{\gamma-2\alpha^{2}}\leq1-\frac{2\alpha}{3}.
\end{align*}
In the above derivation we used the inequalities
\[
\alpha\leq1/4,\text{ \ \ }2K/n^{2}<\alpha^{2},\text{ \ \ }1\geq\gamma
\geq1/2-\alpha/4\geq3\alpha/4>2\alpha^{2}.
\]
Next, Lemma \ref{lev2} implies that%
\[
\frac{\mu\left(  G_{i+1}\right)  }{n-i-1}\geq\frac{\mu\left(  G_{i}\right)
}{n-i}\left(  1+\frac{1}{n-i-1}\left(  \frac{2\alpha}{3}-\frac{1}%
{n-i-1}\right)  \right)  \geq\frac{\mu\left(  G_{i}\right)  }{n-i}\left(
1+\frac{3\alpha}{5n}\right)  .
\]
Therefore,%
\[
\frac{\mu\left(  G_{i+1}\right)  }{n-i-1}\geq\left(  1+\frac{3\alpha}%
{5n}\right)  \left(  1+\frac{3i\alpha}{5n}\right)  \frac{\mu\left(  G\right)
}{n}\geq\left(  1+\frac{3\left(  i+1\right)  \alpha}{5n}\right)  \frac
{\mu\left(  G\right)  }{n},
\]
completing the induction step and the proof of (\ref{in}).

Inequality (\ref{in}) implies that
\begin{align*}
\frac{\mu\left(  \left\vert H\right\vert \right)  }{\left\vert H\right\vert }
&  =\frac{\mu\left(  G_{k}\right)  }{n-k}\geq\left(  1+\frac{3k\alpha}%
{5n}\right)  \frac{\mu\left(  G\right)  }{n}\geq\left(  1+\frac{3k\alpha}%
{5n}\right)  \left(  \gamma-\frac{K}{n^{2}}\right) \\
&  =\gamma\left(  1+\frac{3k\alpha}{5n}\right)  -\frac{K}{n^{2}}\left(
1+\frac{3k\alpha}{5n}\right)  >\gamma\left(  1+\frac{4k\alpha}{7n}\right)
+\frac{\alpha}{42n}-\frac{2K}{n^{2}}\\
&  >\gamma\left(  1+\frac{4k\alpha}{7n}\right)  ,
\end{align*}
as claimed.

To complete the proof of the theorem, note that, after the procedure
$\mathcal{P}$ stops, we have either $k=\left\lfloor \beta n\right\rfloor $ or
$\delta\left(  H\right)  >\left(  \gamma-\alpha\right)  \left\vert
H\right\vert .$ If $k=\left\lfloor \beta n\right\rfloor ,$ then
\[
\mu\left(  H\right)  \geq\gamma\left(  1+\frac{4\left\lfloor \beta
n\right\rfloor \alpha}{7n}\right)  \left\vert H\right\vert >\gamma\left(
1+\frac{\beta\alpha}{2}\right)  \left\vert H\right\vert ;
\]
hence, condition \emph{(i)} holds.

If $k<\left\lfloor \beta n\right\rfloor ,$ then $\delta\left(  H\right)
>\left(  \gamma-\alpha\right)  \left\vert H\right\vert ,$ and, in view of
(\ref{bnd}), we find that%
\[
\mu\left(  H\right)  >\gamma\left(  1+\frac{k\alpha}{2n}\right)  \left\vert
H\right\vert >\gamma\left\vert H\right\vert ;
\]
hence, condition \emph{(ii)} holds, completing the proof.
\end{proof}

\subsubsection*{Proof of Theorem \ref{th1}}

\begin{proof}
[ ]Let $G$ be a graph of order $n$ with $\mu\left(  G\right)  >\sqrt
{\left\lfloor n^{2}/4\right\rfloor }.$ Assume first that $\delta\left(
G\right)  >2n/5.$ Since $G$ contains a triangle, it is nonbipartite; hence,
for $n$ sufficiently large, Fact \ref{tNS} implies that $C_{t}\subset G$ for
every $t\leq\delta\left(  G\right)  +1,$ completing the proof.

Thus, we shall assume that $\delta\left(  G\right)  \leq2n/5.$ Let
\[
\alpha=1/10,\text{ \ \ }\beta=1/2,\text{ \ \ }\gamma=1/2,\text{ \ \ }K=1.
\]
We have $\delta\left(  G\right)  \leq\left(  \gamma-\alpha\right)  n$ and
\[
\mu\left(  G\right)  \geq\sqrt{\left\lfloor n^{2}/4\right\rfloor }\geq
n/2-1/n=\gamma n-K/n.
\]
Hence, Theorem \ref{thv3} implies that, for $n$ sufficiently large, there
exists an induced subgraph $H\subset G$ with $\left\vert H\right\vert >n/2,$
satisfying one of the following conditions:

\emph{(i)} $\mu\left(  H\right)  >\left(  1/2+1/80\right)  \left\vert
H\right\vert ;$

\emph{(ii)} $\mu\left(  H\right)  >\left\vert H\right\vert /2$ and
$\delta\left(  H\right)  >2\left\vert H\right\vert /5.$

Assume first that condition \emph{(i)} holds. Then, by Fact \ref{leNSMM}, we
obtain
\[
k_{3}\left(  H\right)  >\left(  \frac{\mu\left(  H\right)  }{\left\vert
H\right\vert }-\frac{1}{2}\right)  \frac{1}{12}\left\vert H\right\vert
^{3}\geq\frac{1}{80\cdot12}\left\vert H\right\vert ^{3}=\frac{1}%
{960}\left\vert H\right\vert ^{3}.
\]
Thus, there is a vertex $u\in V\left(  H\right)  $ contained in at least
$3k_{3}\left(  H\right)  /\left\vert H\right\vert \geq\left\vert H\right\vert
^{2}/320$ triangles in $H$, and so the neighborhood of $u$ induces more than
$\left\vert H\right\vert ^{2}/320$ edges. By a theorem of Erd\H{o}s and Gallai
\cite{ErGa59}, the neighborhood of $u$ contains a path $P$ longer than
\[
\frac{2}{320}\left\vert H\right\vert \geq\frac{1}{320}n.
\]
Clearly, the path $P$ and the vertex $u$ form a cycle $C_{t}$ for every $t\leq
n/320,$ completing the proof in this case.

If condition \emph{(ii)} holds then, by $\mu\left(  H\right)  >\left\vert
H\right\vert /2,$ the graph $H$ contains a triangle; thus, by Fact \ref{tNS},
$C_{t}\subset H$ for every $t\leq\delta\left(  H\right)  +1,$ completing the proof.
\end{proof}

\subsubsection*{Proof of Theorem \ref{th2}}

\begin{proof}
[ ]Let $G$ be a graph of order $n$ with $\mu\left(  G\right)  >\left(
1/2-\theta\right)  n.$ If $G$ is triangle-free, the proof is completed by Fact
\ref{tstab}, so we shall assume that $G$ contains a triangle.

Assume first that $\delta\left(  G\right)  >2n/5.$ Since $G$ is nonbipartite,
for $n$ sufficiently large, Fact \ref{tNS} implies that $C_{t}\subset G$ for
every $t\leq\delta\left(  G\right)  +1,$ completing the proof.

Thus, we shall assume that $\delta\left(  G\right)  \leq2n/5.$ Let
\[
\alpha=1/10+\theta,\text{ \ \ }\beta=40\theta,\text{ \ \ }\gamma
=1/2-\theta,\text{ \ \ }K=0.
\]
We have $\delta\left(  G\right)  \leq\left(  \gamma-\alpha\right)  n$ and
$\mu\left(  G\right)  >1/2-\theta=\gamma n.$ Hence, Theorem \ref{thv3} implies
that, for $n$ sufficiently large, there exists an induced subgraph $H\subset
G$ with $\left\vert H\right\vert >\left(  1-\beta\right)  n,$ satisfying one
of the following conditions:

\emph{(i)} $\mu\left(  H\right)  >\gamma\left(  1+\alpha\beta/2\right)
\left\vert H\right\vert ;$

\emph{(ii)} $\mu\left(  H\right)  >\gamma\left\vert H\right\vert $ and
$\delta\left(  H\right)  >\left(  \gamma-\alpha\right)  \left\vert
H\right\vert .$

Assume first that condition \emph{(i)} holds. Then,
\begin{align*}
\mu\left(  H\right)   &  \geq\gamma\left(  1+\alpha\beta/2\right)  \left\vert
H\right\vert =\left(  1/2-\theta+\left(  1/2-\theta\right)  \left(
1/10+\theta\right)  20\theta\right)  \left\vert H\right\vert \\
&  =\left(  1/2-\theta-\theta+8\theta^{2}-20\theta^{3}\right)  \left\vert
H\right\vert >\left\vert H\right\vert /2,
\end{align*}
and so, by Theorem \ref{th1}, $C_{t}\subset H\subset G$ for every
$t<\left\vert H\right\vert /320.$ This completes the proof in view of%
\[
\left\vert H\right\vert /320\geq\left(  1-40\theta\right)  n/320>n/321.
\]

If condition \emph{(ii)} holds then, in view of $\delta\left(  H\right)
>\left(  \gamma-\alpha\right)  \left\vert H\right\vert =2\left\vert
H\right\vert /5,$ Fact \ref{tNS} implies that $C_{t}\subset H$ for all
$t\leq\delta\left(  H\right)  +1,$ unless $H$ is bipartite. To complete the
proof we have to consider case of bipartite $H$. Since $H$ is triangle-free
and $\mu\left(  H\right)  >\gamma\left\vert H\right\vert =\left(
1/2-\theta\right)  \left\vert H\right\vert ,$ Fact \ref{tstab} implies that
there exists an indiced bipartite subgraph $G_{0}\subset H$ satisfying
\[
\left\vert G_{0}\right\vert >\left(  1-3\theta^{1/3}\right)  \left\vert
H\right\vert \geq\left(  1-3\theta^{1/3}\right)  \left(  1-\beta\right)
n=\left(  1-3\theta^{1/3}\right)  \left(  1-40\theta\right)  n>\left(
1-4\theta^{1/3}\right)  n
\]
and
\[
\delta\left(  G_{0}\right)  >\left(  1-6\theta^{1/3}\right)  \left\vert
H\right\vert \geq\left(  1-6\theta^{1/3}\right)  \left(  1-\beta\right)
n=\left(  1-3\theta^{1/3}\right)  \left(  1-40\theta\right)  n>\left(
1-7\theta^{1/3}\right)  n,
\]
completing the proof.
\end{proof}

\section*{\textbf{Concluding remarks}}

It is clear that the constant $1/320$ in Theorem \ref{th1} can be increased
even with the present methods; thus, the following question arises:\bigskip

\textbf{Question }\emph{What is the maximum }$C$\emph{ such that for all
positive }$\varepsilon<C$\emph{ and sufficiently large }$n,$\emph{ every graph
}$G$\emph{ of order }$n$\emph{ with }$\mu\left(  G\right)  >\sqrt{\left\lfloor
n^{2}/4\right\rfloor }$\emph{ contains a cycle of length }$t$\emph{ for every
}$t\leq\left(  C-\varepsilon\right)  n.$\bigskip

It is known (\cite{Bol78}, p. 150) that if $G$ is a graph of order $n$ with
$e\left(  G\right)  >\left\lfloor n^{2}/4\right\rfloor $, then $G$ contains a
cycle of length $t$ for every $3\leq t\leq\left\lceil n/2\right\rceil .$ Thus,
one can conjecture that $C=1/2$. However, this is not true: taking the join of
a complete graph of order $k=\left\lceil \left(  3-\sqrt{5}\right)
n/4\right\rceil $ and an empty graph of order $n-k,$ we obtain a graph $H$ of
order $n$ with $\mu\left(  H\right)  >n/2\geq\sqrt{\left\lfloor n^{2}%
/4\right\rfloor },$ but having no cycles longer than $2k\sim\left(  3-\sqrt
{5}\right)  n/2.$\bigskip

Finally, a word about the project mentioned in the introduction: in this
project we try to follow the following principles:

- give results that can be used as wide-range tools, like Lemmas \ref{lev1}
and \ref{lev2}, Theorem \ref{thv3}, and Facts \ref{leNSMM} and \ref{tstab};

- give explicit conditions for the parameters in statements, like the
conditions for $\alpha,\beta,\gamma,K,n$ in Theorem \ref{thv3};

- prefer simple to optimal bounds, like the factor $1/320$ in Theorem
\ref{th1}.

We aim to give results that can be used further, hoping to add more integrity
to spectral extremal graph theory.

\end{document}